\documentclass[11pt,letterpaper]{amsart}

\usepackage[english]{babel}

\usepackage{amsfonts, amssymb, amsmath, eucal, verbatim, amsthm, amscd, enumerate}

\usepackage[a4paper,margin=3cm,innermargin=3cm]{geometry}

\usepackage[colorlinks=true, allcolors=blue]{hyperref}

\newtheorem{theorem}{Theorem}

\newtheorem{proposition}[theorem]{Proposition}

\newtheorem{fact}[theorem]{Fact}

\theoremstyle{definition}

\def\R{{\mathbb R}}

\renewcommand{\R}{\mathbb{R}}

\begin{document}

\title{Duality and Heat flow}

\author{D. Cordero-Erausquin}
\address{D. C.-E.: Institut de Mathématiques de Jussieu, Sorbonne Université, CNRS, F-75005 Paris. }

\author{N. Gozlan}
\thanks{N.G is supported by a grant of the Agence nationale de la recherche (ANR), Grant ANR-23-CE40-0017 (Project SOCOT)}
\address{N. G.:  Université Paris Cité, CNRS, MAP5, F-75006 Paris.}

\author{S. Nakamura}
\thanks{S. N. was supported by JSPS Overseas Research Fellowship and JSPS Kakenhi grant number 21K13806} 
\address{S. N.: Department of Mathematics, Osaka University, Osaka, Japan.}

\author{H. Tsuji}
\thanks{H. T. was partially supported by JSPS Kakenhi grant number 22J10002}
\address{H. T.: Department of Mathematics, Saitama University, Saitama, Japan.}

\begin{abstract}
    We reveal the relation between the Legendre transform of convex functions and Heat flow evolution, and how it applies to the functional Blaschke-Santal\'{o} inequality. We also describe local maximizers in this inequality. 
\end{abstract}

\maketitle

\section{Introduction}
Our goal in this note is twofold. Firstly, we will
describe how duality for log-concave functions commutes with Gaussian semi-groups, like Heat flow. Secondly, we will use this commutation property to give a simple, streamlined semi-group proof of the functional form of the Blaschke-Santal\'o inequality for centrally symmetric convex bodies, in the wake of the approach proposed by the two last authors. This will allow us to describe local maximizers in the inequality. 


Before we elaborate more on these points, let us recall  the functional form of the Blaschke-Santal\'o inequality: for all even (convex) function $\phi:\R^n\to \R\cup\{+\infty\}$, 
\begin{equation}\label{eq:funct_s1}
M(e^{-\phi}):=\int_{} e^{-\phi} \int_{} e^{-\phi^\ast} \le (2\pi)^{n}. 
\end{equation}  
The left-hand side term is called the functional volume product and the inequality states it is maximized when $\phi=|\,\cdot\,|^2/2$, the fixed point of the Legendre's transform, $$\phi^\ast(z):=\sup_{x} x\cdot z - \phi(x) \in \R\cup\{+\infty\}, \qquad \forall z\in \R^n,$$
on the Euclidean space $(\R^n, |\,\cdot\,|)$. 
Taking $\phi=\|\,\cdot\,\|_K^2 /2$, where $\|\,\cdot\,\|_K$ is the gauge associated to a centrally symmetric convex body $K\subset \R^n$,  
the geometric inequality 
\[
{\rm vol}(K) {\rm vol}(K^\circ) \le {\rm vol}(B_2^n)^2
\]
is recovered; here $K^\circ = \{x\in \R^n\; ; \; x\cdot y \le 1, \; \forall y \in K\}$ and $B_2^n =\{|\,\cdot\,|\le 1\}$ (see e.g.~\cite{MP, Lut}).

Inequality~\eqref{eq:funct_s1} was put forward by Keith Ball in his PhD dissertation. It has been extended to non-even functions in~\cite{AKM} and admits several proofs relying on geometric arguments (see for instance~\cite{AKM, FM:2007, L:2009}). 
Let us mention that one may prefer to work directly with log-concave functions, that is functions of the form $e^{-\phi}$ with $\phi:\R^n\to \R\cup\{+\infty\}$ convex. Then, the \emph{polar} of the log-concave function $f=e^{-\phi}$ is defined by
$$f^\circ(z):= e^{-\phi^\ast(z)}=\inf_{x\in \R^n}\frac{e^{-x\cdot z}}{f(x)},\qquad \forall z\in \R^n$$
so, among even log-concave functions, $M(f)=\int_{}  f \, \int_{}  f^\circ$ is maximized when $f$ is a centered Gaussian, that is a function of the form $f(x)=C\, e^{-Hx\cdot x}$ where $H$ is a definite positive matrix and $C>0$. Recall $M$ is linearly invariant, in the sense that if $\tilde f= \lambda f\circ T$, where $T\in GL_n(\R)$ and $\lambda>0$, then $M(\tilde f)= M(f)$.

A new, analytical, approach to~\eqref{eq:funct_s1} based on semi-group interpolation was obtained in~\cite{NT:2023}, where it is proven  that $M(e^{-\phi})$ increases when we let $e^{-\phi}$ evolve along the Fokker-Planck flow. 
The idea in there is to regard the polar transform as the limit of a suitably rescaled Laplace transform. That is, 
\[
\lim_{p\to0^+} \big( \int_{} e^{\frac{1}p x\cdot z} f(x)^{\frac{1}{p}}\, dx \big)^{-p}
= 
f^\circ(z), \qquad \forall z \in \R^n,
\]
see \cite{NT:2023, CFL:2024} for more details. 
Thus, the study of the polar transform may be reduced to that of some '$p$-Laplace transform', $p\in (0,1)$. %
Based on this idea, the last two authors first established the monotonicity of some variant of the volume product along Heat flow (with a $p$-Laplace transform in place of the polar transform) by use of the variance Brascamp-Lieb inequality, recalled below, and the Cram\'{e}r-Rao inequality. 
The monotonicity of the volume product along Heat flow was then obtained in the limit as $p \to 0^+$.
In contrast, our approach to the monotonicity does not require such limiting argument. Our key tools will be a new formula for the evolution of the polar of a log-concave function, combined with the Brascamp-Lieb inequality, only. Thus our work makes the approach to the Blaschke-Santal\'o inequality of the prior work clearer, both at technical and conceptual levels. Moreover, it will allow us to prove that Gaussian functions are the only local maximizers, a result that cannot be reached with a limiting argument.

\section{Heat flow and Legendre transform}

Let us first recall some simple and well known properties of log-concave functions (see e.g. \cite[Section 2]{CFL:2024} and references therein).  We shall be considering log-concave functions $f=e^{-\phi}$ with $0<\int f <\infty$. Note for consistency that, for an even log-concave function $f$ we have
\begin{equation}\label{eq:finite}
    \int f >0 \Rightarrow\int f^\circ <\infty.
\end{equation}
The argument is standard: if we denote by $C$ the interior of the set $\{f>0\}$, then $C$ is a convex set, where $f$ is continuous, and by assumption $C$ must have non-empty interior and thus contains a small ball $B(r)$ around $0$; this implies that $f\ge M 1_{B(r)}$ for some constant $M>0$, which in turn gives $f^\circ (x) \le M^{-1} \, e^{-|x|/r}$, $x \in \R^n.$

We say that $\phi$ is \emph{coercive} if  $\phi(x)\ge  a|x| - b$ for all $x\in \R^n$ for some constants $a>0$, $b\in \R$. This assumption is natural here since for a convex function $\phi$,
$$\phi \textrm{ coercive} \Longleftrightarrow \int e^{-\phi} < \infty.$$
\noindent To deal with duality, we shall assume that $\phi$ is \emph{super-linear} which means that $\limsup_\infty \frac{\phi(x)}{|x|}=\infty$. This is equivalent to the fact the domain of $\phi^\ast$ is $\R^n$, that is 
\begin{equation} \label{eq:superlinear}
\phi \textrm{ super-linear} \Longleftrightarrow
\phi^\ast \textrm{ takes finite values only} \Longleftrightarrow e^{-\phi^\ast}>0,
\end{equation}
see e.g \cite[Propositions 1.3.8 and 1.3.9]{HUL}.
Extending our result below, and understanding the evolution of the domain, for non-super-linear convex functions is a challenging technical question that we will not address here. 

Given a Borel nonnegative function $u$ on $\R^n$, its  Heat flow evolution $P_t u$ 
 is given for $t>0$  by
 \begin{equation}\label{eq:defHeat}
P_t u (x) =\int u(y) \, e^{- |x-y|^2/4t} \frac{dy}{(4\pi t)^{n/2}} = u\ast \gamma_t (x) = \int u(x+\sqrt{2t} y)\, \gamma(y)\, dy, \quad \forall x\in \R^n
 \end{equation}
where $\gamma_t (x)=(\sqrt{2t})^{-n} \gamma(x/{\sqrt{2t}})$, $x\in \R^n$, is an approximate identity given by the Gaussian density $\gamma(x)= (2\pi)^{-n/2}\, e^{-|x|^2/2}$, $x\in \R^n.$
If $u\in L^1(\R^n)$ this defines a $C^\infty$ function $(t,x)\to P_tu (x)$  on $(0,\infty)\times \R^n$, with  $\partial_t (P_t u)(x) = \Delta (P_t u)(x)$. As $t\to 0$, we have convergence of $P_t u$  in $L^1(\R^n)$ to $u$, and if we further assume $u$ bounded  we have pointwise convergence at every point where $u$ is continuous.  When $u$ is an integrable log-concave function, we thus have $P_t u (x) \to u(x)$ at almost-every $x\in \R^n.$

It is often convenient to rescale the Heat evolution in order to obtain an asymptotic profile. 
The Fokker-Planck evolution $Q_t u$ is for $(t,x)\in (0,+\infty)\times \R^n$ given by
\begin{equation}\label{eq:defFP}
Q_t u(x) = e^{nt} P_{\frac{e^{2t}-1}{2}} u(e^{t}\, x)
=\left(\int u(y) e^{\frac{e^{t}}{e^{2t}-1} x\cdot y - \frac1{2(e^{2t}-1)}|y|^2} \, dy\right)\, \frac{e^{-\frac1{1-e^{-2t}} |x|^2/2}}{(2\pi(1-e^{-2t}))^\frac{n}{2}} . 
\end{equation}
It therefore satisfies
$\partial_t (Q_t u) (x) = \Delta Q_t u(x) + {\rm div} (x\, Q_t u)(x)$.  Throughout the paper, the derivatives $\nabla, D^2, {\rm div}, \Delta$ refer to derivatives in the space variables in $\R^n$.

If $f$ is log-concave, then it is well known (for instance by Prékopa's theorem \cite{BL:1976}, or the earlier work~\cite{DKH})  that, for any fixed $t>0$, $P_t f$ is log-concave, strictly positive when $\int f>0$. Actually, we will include for completeness a proof in the Appendix of the following fact.

\begin{fact}\label{fact}
Let $f=e^{-\phi}$ a log-concave function with $0<\int f <\infty$, and for $t>0$ fixed, denote $\phi_t := -\log P_t f$. Then 
$D^2 \phi_t(x) >0$ at every $x\in \R^n$, and moreover $\phi_t$ is super-linear when $\phi$ is super-linear.  
\end{fact}

With the notation of Fact \ref{fact}, we see that the evolution of $\phi_t$ is given by 
\begin{equation}\label{eq:logHeat}
\partial_t \phi_t(x) = \Delta \phi_t (x) - |\nabla \phi_t(x)|^2.
\end{equation}
Our main observation is to describe the evolution of $(\phi_t)^\ast$, which is simple  but most useful, as we shall see.

\begin{proposition}\label{prop:main}
Let $\phi:\R^n \to \R\cup\{+\infty\}$ be a super-linear convex function with $\int e^{-\phi}>0$. Consider the convex function $\phi_t = - \log P_t (e^{-\phi})$ where $P_t$  is the  Heat semi-group and let $\psi_t = (\phi_t)^*$ be its Legendre transform. Then for every $z\in \R^n$ and $t>0$
\begin{equation}\label{e:PointwiseId}
    \partial_t \psi_t  (z) = |z|^2 - {\rm Tr }(D^2\psi_t(z))^{-1} . 
\end{equation}
\end{proposition}

\begin{proof}
The function $\phi_t$ is smooth, super-linear and strictly convex, from Fact \ref{fact} above.  In particular,  $\psi_t=(\phi_t)^\ast$ is finite and smooth on $\R^n$, and $\nabla\phi_t$ is a diffeomorphism of $\R^n$ with inverse $\nabla \psi_t$.  We next use the following simple relation, valid for any first order perturbation of a convex function:
$$\partial_t \psi_t(z) = - \partial_t \phi_t (\nabla \psi_t (z)).$$
To check this, take for instance the derivative in $t$ of $\psi_t (z) + \phi_t(\nabla \psi_t (z)) = z\cdot\nabla\psi_t(z)$ and use that $\nabla\phi_t (\nabla\psi_t(z))=z$. In the case of Heat flow~\eqref{eq:logHeat} we have
$$\partial_t \phi_t(y) = \Delta \phi_t(y)  - |\nabla \phi_t(y) |^2 = {\rm Tr}D^2 \phi_t(y)  - |\nabla \phi_t(y) |^2 .$$
The relation~\eqref{e:PointwiseId} follows. 
\end{proof}

By elementary change of variables or change of functions, one can state  similar formulas in the case of the Fokker-Planck flow, or for the Gaussian reformulation in terms of infimal convolution $f^c(y):=-(|x|^2/2 + f)^\ast(y)-|y|^2/2$  (see~\cite{M:1991, BGL} for background) when $e^{-f}$ evolves along the Ornstein-Uhlenbeck semi-group, which means that $e^{-|x|^2/2-f(x)}$ follows the Fokker-Planck semi-group.

For instance, if we use the Fokker--Planck evolution $\phi_t:= -\log\, Q_tf$ instead of the Heat evolution, we have 
$
\partial_t \phi_t(x) = \Delta \phi_t (x) - |\nabla \phi_t|^2(x) + x\cdot \nabla \phi_t (x) - n, 
$
and, either by the above argument or by change of variables, we have for  $\psi_t:= (\phi_t)^*$ that
$$
\partial_t \psi_t(z)
= 
|z|^2 - {\rm Tr}(D^2 \psi_t(z))^{-1} - z\cdot \nabla \psi_t(z) +n. 
$$

\section{Blashke-Santal\'o inequality and Heat flow}

We now use Proposition~\ref{prop:main} to provide a simple proof of the monotonicity of the volume product $M$ along the  Fokker-Planck semi-group; this in turn proves ~\eqref{eq:funct_s1}. As explained in the Introduction, such monotonicity was obtained in~\cite{NT:2023}, without even the assumption of log-concavity, actually, as a limit case of the monotonicity of another functional involving the Laplace transform instead of the polar transform. Our direct method will allow us to give a condition for strict-monotonicity, and then to describe all local maximizers of $M$.

\begin{theorem}\label{theo1}
Let $f=e^{-\phi}$ be an even log-concave function with $0<\int f<\infty$. Then, if $f_t$ denotes either the Heat or Fokker-Planck evolution of $f$, we have that
$$\alpha(t)=M(f_t)$$
increases in $t\in[0,+\infty)$, and, assuming $\phi$ is super-linear, it  strictly increases unless $f$ is a centered Gaussian function. 
\end{theorem}

\begin{proof}
By the linear invariance of $M$, we see that
$$M(Q_tf) = M(P_{\ell(t)} f),\qquad \forall t\geq0,$$
for some strictly increasing function $\ell(t)$, so we can work with the Heat semi-group only  (this reduction was noted in~\cite{CFL:2024}). Let $\phi_t = -\log P_t f$, i.e. $f_t = e^{-\phi_t}$, and $\psi_t = (\phi_t)^\ast$, so that 
$$\alpha(t):=\log M(P_tf)=\log\int e^{-\phi}  + \log \int e^{-\psi_t}.$$

We fist prove that $\alpha$ increases.
The heuristic for monotonicity is as follows: we formally have, using the previous Proposition, that
$$\alpha'(t) = \int  (  {\rm Tr} (D^2_x\psi_t)^{-1}  - |x|^2 ) \frac{e^{-\psi_t(x)} dx}{\int e^{-\psi_t} },$$
and so summing for $i=1,\ldots n$, the Brascamp-Lieb inequality~\eqref{eq:BL} below applied to the linear functions $u(x)=x_i$, for which $\int u \, e^{-\psi_t}=0$ since $\psi_t$ is even, we arrive at  $\alpha'(t)\ge 0$. 

Let us now give the rigorous arguments. We are given a log-concave even function $f=e^{-\phi}$ with $\phi$ coercive and $\int f>0$. 
Let us first note that it suffices to consider super-linear functions $\phi$. Indeed, if we introduce the log-concave functions defined, for all $x\in \R^n$, by 
$$f_k(x) = f(x) e^{-|x|^2/k} = e^{-\phi(x) - |x|^2/k} \le f (x),$$
then the function $-\log f_k$ is super-linear. By dominated  convergence, we have, for fixed $t>0$, $F_k:=P_t f_k \nearrow P_t f=:F$ pointwise, as $k\to \infty$. Thus, by monotone convergence, $\int F_k \to \int F$. On the other hand, the sequence $(F_k)^\circ$ decreases, and $\int (F_k)^\circ <\infty$ since $\int F_k >0$. Thus, by dominated convergence,
$\int (F_k)^\circ \to \int \inf_k (F_k)^\circ$. But
$$ \inf_k (F_k)^\circ (x) = \inf_k \inf_y \frac{e^{-x\cdot y}}{F_k(y)} =  \inf_y \frac{e^{-x\cdot y}}{\sup_k F_k(y)} = F^\circ (x),\qquad \forall x\in \R^n.$$
This shows that $M(P_t f_k) \to M(P_t f)$, allowing to reduce monotonicity to the case the convex functions are super-linear. 

So, in the sequel, we assume $\phi$ is an even super-linear convex function with $\int e^{-\phi}>0$. 
It is easily checked that $\partial_t e^{-\psi_t}$ is locally uniformly in $t>0$ dominated by an integrable function,  by virtue of Proposition \ref{prop:main}; indeed, the same arguments giving~\eqref{eq:hessian} give that $D^2 \phi_t \le \frac1{2t}$ (see also~\cite[Lemma 1.3]{EL}) and so $D^2 \psi_t \ge 2t$.
So we have, using Proposition \ref{prop:main}, that 
$$\alpha'(t) = \int ( -\partial_t \psi_t)\,   \frac{e^{-\psi_t}}{\int e^{-\psi_t} }= \int  (  {\rm Tr} (D^2\psi_t(x))^{-1}  - |x|^2 ) \frac{e^{-\psi_t(x)} dx}{\int e^{-\psi_t} },\qquad \forall t>0.$$
We now call upon  the variance Brascamp-Lieb inequality~\cite{BL:1976} which states that for a $C^2$-smooth convex function $V:\R^n \to \R$ with $D^2 V>0$
almost-everywhere, denoting by $\mu_V$ the probability on $\R^n$ with density $\frac{e^{-V}}{\int e^{-V}}$, the inequality
\begin{equation}\label{eq:BL}
\textrm{Var}_{\mu_V} (u)= \int \Big( u - \int u\, d\mu_V\Big)^2 \, d\mu_V \le \int (D^2 V)^{-1} \nabla u \cdot \nabla u \, d\mu_V,
\end{equation}
holds for every smooth $u\in L^2(\mu_V)$. 
We apply this inequality with $V=\psi_t$ and 
to the linear functions $u:x\to x_i$, $i\le n$, which are centered, $\int x_i e^{-\psi_t(x)}\, dx =0$,  since $\psi_t$ is even. Summing over $i\le n$, we then get exactly that $\alpha'(t)\ge 0$, for $t>0$. 
At $t=0^+$ we have
\begin{equation}\label{eq:lim0}
\liminf_{t\to 0^+} \int e^{-\psi_{t}} \ge \int e^{-\psi}. 
\end{equation}
Indeed, for every $x,y,z\in \R^n$ we have
$-\phi(y) - |y-x|^2/4t \le \psi(z) - z\cdot y - |y-x|^2/4t$ and therefore
$$ \phi_t(x) \ge -\psi(z) + z\cdot x -  t |z|^2.$$
This implies, by definition, that for every $z\in \R^n$, 
$\psi_t(z) \le \psi(z) + t |z|^2$,
and so 
$
\int e^{-\psi_{t}}\ge
\int e^{-\psi(x) - t x^2}\, dx.$
Then~\eqref{eq:lim0} follows by Fatou's Lemma. 
So we have proved that $\alpha$ increases on $[0,+\infty)$. 

Finally, assume that $\alpha$ is not strictly increasing, which implies that $\alpha'(t)=0$ for some $t>0$. We see that in the argument above (for this fixed $t$), we must have  equality in the Brascamp-Lieb inequality for all linear functions $x\to x_j$. However, the only centered equality cases in~\eqref{eq:BL} are the linear combinations of the $\partial_i V$, $i\le n$, which forms a space of dimension at most $n$; this is a classical fact (one can call upon the discussion in the last section of~\cite{CE:2017} for instance).  Since the linear functions $x\to x_j$ are linearly independent, it means that each $\partial_i \psi_{t} $ is a linear combination of the $x\to x_j$, $j\le n$. This implies that $D^2\psi_{t}$ is constant on $\R^n$. Since $\psi_{t}$ is even, it means that $\psi_{t}$ and thus $\phi_{t}=\psi_{t}^\ast$ are quadratic functions (up to constants), so $\phi_t(x)=\frac12 Hx\cdot x + C$ with $H$ a positive definite matrix and $C\in \R$. 
We conclude with a standard analytic argument: if $f_t=e^{-\phi_{t}} = P_{t} f$ is a centered Gaussian function, then $f$ must be so too. Indeed, if $\hat g(\xi) = \int g(x) e^{-ix\cdot \xi}\, dx$ denotes the Fourier transform of $g\in L^1(\R^n)$, we have that $\widehat{f_t}(\xi) = \hat f (\xi) \, e^{-t |\xi|^2}$. Since $\widehat{f_t}(\xi) = \tilde C e^{-H^{-1} \xi \cdot \xi}$ for some constant $\tilde C>0$, we have that $\hat f(\xi) = \tilde C \, e^{-A\xi \cdot \xi}$ for some symmetric matrix $A$. But $\hat f$ has to tend to $0$ at infinity, since $f\in L^1(\R^n)$, thus $A$ must be positive definite. By injectivity of the Fourier transform on $L^1$, we deduce that $f$ is a centered Gaussian function. 
\end{proof}

In order to recover the Blaschke-Santal\'o inequality, it suffices to note that, with the notation of Theorem \ref{theo1}, when using the Fokker-Plank flow, we have $M(Q_t f)\to M(\gamma)$ as $t\to \infty$. Indeed, using the bound $f(y)\le C e^{-c|y|}$, we see from~\eqref{eq:defFP} and dominated convergence that $Q_t(f)(x) \to \big(\int f\big) \gamma(x)$ at every fixed $x$.  By~\cite[Lemma 3.2]{AKM}, it  implies that  $ (Q_t(f))^\circ (y) \to \big(\int f\big)^{-1}(2\pi)^n \gamma(y)$ at every $y$, since $\gamma^\circ=(2\pi)^n\gamma>0$, and these pointwise convergences imply in turn the convergence of the integrals (for this, see also~\cite[Fact 2.5]{CFL:2024}). 

The strict-monotonicity will allow us to prove the following result, which can be viewed as a functional version of the geometric result in~\cite{MR}.

\begin{theorem}
The functional $f\to M(f)$ on the subset of $L^1(\R^n)$ formed by even log-concave functions  has no local maximizers besides centered Gaussian functions, which are the global maximizers. 
\end{theorem}

\begin{proof}
If the local maximizer is of the form $e^{-\psi}$ with $\psi$ super-linear, we can invoke the strict monotonicity from the previous Theorem to conclude. There is a trick that allows to use this strict-monotonicity in general, as we now explain. 

So assume that $f$ is an even log-concave function that is local maximizer of $M$: for some $\epsilon>0$, 
\vskip-0.5cm
$$M(f) = \sup\Big\{ M(g)\; ; \; g \textrm{  even log-concave with }
\|f-g\|_{L^1(\R^n)} \le \epsilon \Big\}.$$
First note that $0<\int f < \infty$, since $\int f = \infty \Rightarrow \int f^\circ=0$ by~\eqref{eq:finite} applied to $f^\circ$. Since $P_tf \to f$ in $L^1(\R^n)$ as $t\to 0$, let us fix a small $t_0>0$, so that
$g:= P_{t_0} f$ verifies $\| f - g \| \le \epsilon/2$.
Since $M(g)\ge M(f)$ by the previous Theorem, we must have $M(g)=M(f)$. But note that $g$ is continuous and strictly positive on $\R^n$, therefore, if we consider the even log-concave function $G:=g^\circ$, we have $G= e^{-\psi}$ for some convex super-linear function $\psi$, in view of~\eqref{eq:superlinear}, and $M(G)=M(f)>0$. 

We now let $G$ evolve along the Heat semi-group. For $t$ small $(P_t G)^\circ$ will remain at $L^1$ distance from $f$ smaller than $\epsilon$, because $(P_t G)^\circ \to g$ in $L^1(\R^n)$ as $t\to 0$. Indeed,  we have that  $P_t G$ converges almost everywhere to $G$. Since  $G^\circ=g>0$, this implies, by~\cite[Lemma 3.2]{AKM},  that $(P_t G)^\circ$ converges to $G^\circ=g$ pointwise on $\R^n$ as $t\to 0^+$, which in turn implies  $\int (P_t G)^\circ \to \int G^\circ $ (see~\cite[Fact 2.5]{CFL:2024}). But, according to Scheffé's lemma, pointwise convergence and convergence of integrals, imply, for nonnegative functions, convergence in $L^1(\R^n)$, as wanted.

Thus,  we must have $M(P_tG ) = M((P_t G)^\circ) = M(g)=M(G)$,  for small $t$'s. But since $G$ is super-linear, it means that $G$ must be a centered Gaussian function, and so must $g=G^\circ$. Since $g=P_{t_0} f$, this forces $f$ to be a centered Gaussian function, as we explained at the end of the proof of the previous Theorem. 
\end{proof}

\medskip

Let us conclude with a word on linear invariance. The presence of this large class of invariance might have been seen as an obstacle to a semi-group approach of the Blaschke-Santal\'o inequality. However, we luckily have that the key inequality is the variance Brascamp-Lieb inequality (our formula~\eqref{e:PointwiseId} points right to it, in fact), and this inequality possesses indeed a linear invariance, unlike usual Poincar\'e inequalities, say, which depend on the Euclidean structure. In this regard, it is natural to investigate other linear invariant inequalities,  such as the affine isoperimetric-type inequalities in~\cite{Lut, milm}.
We believe the result discussed in this note is the beginning of a promising direction.

\section{Remarks on the non-even case}

We would like to give some hints on how general log-concave functions can be treated, without entering in all technical details. 

Blaschke-Santal\'o inequality cannot hold for general sets or functions without some form of centering. As put forward in~\cite{AKM}, for log-concave functions $f:\R^n\to \R$ with $0<\int f <\infty$, we have again that Gaussian functions maximize the volume product provided we take the infimum over translations, that is
$$\int f \inf_{z\in \R^n} \int (\tau_z f)^\circ \le (2\pi)^n$$
where $\tau_z f(x)= f(x-z)$. This was derived in~\cite{CFL:2024} by observing that
\begin{equation}
t\to  \beta(t) :=\log\inf_{z\in \R^n} \int (\tau_z f_t)^\circ    
\end{equation}
increases in time $t>0$, where $f_t$ denotes either the Heat or Fokker-Plank evolution of $f$. However, the path in~\cite{CFL:2024} was the same as in~\cite{NT:2023}: the result was first proved for a '$p$-Laplace transform' (see also the independent subsequent work~\cite{Mast}) and then derived by a limiting argument as $p\to 0$. Here, we would like to sketch a direct approach, combining the ideas of~\cite{CFL:2024}  and  the Proposition~\ref{prop:main} above. 

Given a super-linear log-concave function $f$ with $\int f>0$, let us introduce
$$R(t,z) = \log \int (\tau_z f_t)^\circ = \log\int e^{-\psi_t} \, e^{-z\cdot x}\, dx =\log L(e^{-\psi_t})(-z)$$
where $f_t$ is the Heat flow evolution of $f$, $\psi_t := (-\log f_t)^\ast$, as before,   and $L$ stands for the Laplace transform of a non-negative function: $L(g)(z):=\int g(x)\, e^{z\cdot x} \, dx$.

We proceed as above: for fixed $z\in \R^n$,
$$\partial_t R(t,z) = - \int \partial_t \psi_t\, d\mu_{V_z} = \int \Big( {\rm Tr}(D^2\psi_t(x)^{-1}) - |x|^2 \Big)\, d\mu_{V_z}(x)$$
where $V_z(x) := \psi_t(x) + z\cdot x$ is a convex potential on $\R^n$ with $D^2 V_z = D^2\psi_t$. We can again call upon the Brascamp-Lieb inequality~\eqref{eq:BL} for the linear functions $u(x)= x_i$, but now we will have a centering term that remains, 
$$\partial_t R(t,z) \ge -\Big|\int x \, d\mu_{V_z}(x) \Big|^2 = -|\nabla R (t,z)|^2;$$
the equality shows the importance of the Laplace transform in the argument. Thus we have the inequality
$$\partial_t R(t,z) + |\nabla  R(t,z)|^2 \ge 0.$$
The function $z\to R(t,z)$ is convex  and admits a unique minimizer $s(t)\in \R^n$, see e.g.~\cite{CFL:2024}. Thus, we have for our function
$\beta(t)=R(t,s(t))$ that
$$\beta'(t) =\partial_t R(t, s(t)) + s'(t)\cdot \nabla R(t,s(t))
\ge -|\nabla R(t,s(t))|^2 + s'(t)\cdot \nabla R(t,s(t)) = 0$$
since $\nabla R(t,s(t))=0$. This shows that $\beta$ increases, as wanted. Of course, more work should be done in order to justify the formal computations above; our intention here is only to explain how our arguments combine with those of~\cite{CFL:2024} (where the interested reader will find ways to handle technical details).

\section*{Appendix: proof of Fact~\ref{fact}}
We recall the notation $\phi_t= - \log P_t(e^{-\phi})$ for $t>0$ fixed, where $\phi$ is a coercive convex function with $\int e^{-\phi}>0$. We know that $\phi_t$ is convex, but one can notice that this also follows from the computation below, which is the approach of Brascamp-Lieb~\cite{BL:1976}.

We fist prove the strict convexity of $\phi_t$, in the form $D^2 \phi_t >0$. 
Since $e^{-\phi_t} = P_{t/2}(P_{t/2}e^{-\phi}) $ and $P_{t/2}e^{-\phi}$ is an integrable, positive, smooth, log-concave function, we can assume that  $\phi:\R^n\to \R$ is a convex coercive twice continuously differentiable function. We can also assume that $t=1/2$,  for notational simplicity. 
For a fixed direction $|\theta|=1$, we readily check from 
$$e^{-\phi_{1/2} (x)} = \int e^{-\phi(y)} e^{-|x-y|^2/2} \frac{dy}{(2\pi)^{n/2}}
$$
that for any $x\in \R^n$,
\begin{equation}\label{eq:hessian}
(D^2 \phi_{1/2}(x)) \theta\cdot \theta = 1- {\rm Var}_{d\mu_x(y)}(y\cdot \theta)
\end{equation}
where the variance is computed with respect to the probability measure
$$d\mu_x (y) = e^{-\phi(y)} e^{-|x-y|^2/2}\,  \frac{dy}{\int e^{-\phi(z)} e^{-|x-z|^2/2} \, dz}.
$$
Assume that, for some fixed $x\in \R^n$, we have $(D^2 \phi_{1/2}(x)) \theta\cdot \theta = 0$. We then have by the Brascamp-Lieb inequality that
$$1= {\rm Var}_{d\mu_x(y)}(y\cdot \theta) \le
\int \left(D^2 \phi+{\rm Id}_n\right)^{-1}\theta \cdot \theta \,d\mu_x
\le 1$$
since $D^2 \phi+{\rm Id}_n\ge {\rm Id}_n$. As the density of $\mu_x$ is continuous and positive on $\R^n$, we must have $\left(D^2\phi+{\rm Id}_n\right)^{-1}\theta \cdot \theta = 1= |\theta|^2$, and thus $(D^2\phi(y))\theta\cdot \theta = 0$, at every $y \in \R^n$; for this recall that for a nonnegative matrix $H$ we have $\left(H+{\rm Id}_n\right)^{-1} \theta\,  -\,  \theta = -\left(H+{\rm Id}_n\right)^{-1} H \theta$ and that $\left(H+{\rm Id}_n\right)^{-1} H $ is  nonnegative.  So for all $y \in \R^n$, the function $r\mapsto \phi(y + r \theta)$ is affine, which contradicts coercivity, for instance. 

We now prove that when $\phi$ is super-linear, so is $\phi_t$. Since $\phi$ is (bounded-below) super-linear, for an arbitrary $M>0$, there exists $b=b(M)$
such that for every $x\in \R^n$
$$\phi(x) \ge M|x| - b.$$
Therefore we have, using  $|x+\sqrt{2t} y|\ge 
|x| - \sqrt{2t}|y|$
 that
$$e^{-\phi_t}(x)
\le e^{b} e^{-M|x|} \int e^{\sqrt{2t}M |y|} \, d\gamma(y) 
\leq  e^{b} e^{-M|x| + t M^2+ \sqrt{2nt}M  },
$$
where $\gamma$ is the standard Gaussian measure,
that is
$\phi_t(x) \ge M|x| - b - t M^2-\sqrt{2nt}M
$.
Therefore there exists a $K=K(M,n,t)$ such that 
for every $|x|\ge K$, 
$$\frac{\phi_t(x)}{|x|} \ge \frac{M}2 .$$
This implies that $\phi_t$ is  superlinear.

\bibliographystyle{alpha}

\end{document}